\documentclass[10pt]{article}
\usepackage[all]{xy}
\usepackage{amsfonts,amsmath,oldgerm,amssymb,amscd}
\newcommand{\ra}{\rightarrow}

\newcommand{\by}[1]{\stackrel{#1}{\ra}}

\newcommand{\ol}{\overline}		
\newcommand{\iso}{\by \sim}

\newtheorem{theorem}{Theorem}[section]
\newtheorem{proposition}[theorem]{Proposition}
\newtheorem{lemma}[theorem]{Lemma}
\newtheorem{definition}[theorem]{Definition}
\newtheorem{corollary}[theorem]{Corollary}
\newtheorem{question}[theorem]{Question}

\newcommand{\ga}{\alpha}		\newcommand{\gb}{\beta}
			
\newcommand{\gf}{\varphi}

	\newcommand{\BF}{\mbox{$\mathbb F$}}

	\newcommand{\BN}{\mbox{$\mathbb N$}}

\newcommand{\mm}{\mbox{$\mathfrak m$}}

\newcommand{\op}{\mbox{$\oplus$}}	
	
\newcommand{\hh}{\mbox{\rm ht\;}}	\newcommand{\Hom}{\mbox{\rm Hom}}
  	
\newcommand{\Um}{\mbox{\rm Um}}		
\newcommand{\GL}{\mbox{\rm GL}}		
\newcommand{\Aut}{\mbox{\rm Aut}}   
\newcommand{\EL}{\mbox{\rm EL}}

\begin{document} 
\begin{center}
{\large \bf Projective modules over overrings of polynomial
        rings}\\\vspace{.2in} {\large Alpesh M. Dhorajia
        and Manoj K. Keshari}\\
\vspace{.1in}
{\small
Department of Mathematics, IIT Mumbai, Mumbai - 400076, India;\;
        alpesh,keshari@math.iitb.ac.in}
\end{center}

{\small
\noindent{\bf Abstract:} Let $A$ be a commutative Noetherian ring of
dimension $d$ and let $P$ be a projective
$R=A[X_1,\ldots,X_l,Y_1,\ldots,Y_m,\frac {1}{f_1\ldots f_m}]$-module
of rank $r\geq$ max $\{2,\dim A+1\}$, where $f_i\in A[Y_i]$. Then 

$(i)$ $\EL^1(R\op P)$ acts transitively on $\Um(R\op P)$. In
particular, $P$ is cancellative (\ref{23}).

$(ii)$ If $A$ is an affine algebra over a field, then $P$ has a
unimodular element (\ref{222}). 

$(iii)$ The natural map $\Phi_r : \GL_r(R)/EL^1_r(R) \ra K_1(R)$ is
surjective (\ref{201}).

$(iv)$ Assume $f_i$ is a monic polynomial. Then $\Phi_{r+1}$ is an
isomorphism (\ref{201}).

In the case of Laurent polynomial ring (i.e. $f_i=Y_i$), $(i)$ is due
to Lindel \cite{lindel}, $(ii)$ is due to Bhatwadekar, Lindel and
Rao \cite{BLR} and
$(iii,\, iv)$ is due to Suslin \cite{SU}.

\vspace*{.1in}
 \vspace{.1in}

  }

\section{Introduction}
All the rings are assumed to be commutative Noetherian and all the
modules are finitely generated.

Let $A$ be a ring of dimension $d$ and let $P$ be a projective
$A$-module of rank $n$. We say that $P$ is {\it cancellative} if $P\op
A^m\iso Q\op A^m$ for some projective $A$-module $Q$ implies $P\iso
Q$. We say that $P$ has a {\it unimodular element} if $P\iso P'\op A$ for
some projective $A$-module $P'$.

Assume rank $P >\dim A$. Then $(i)$ Bass \cite{Bas} proved that $\EL^1
(A\oplus P)$ acts transitively on $\Um(A\oplus P)$. In particular, $P$
is cancellative and $(ii)$ Serre \cite{Serre} proved that $P$ has a
unimodular element.

Later, Plumstead \cite{P} generalized both the result by proving that
if $P$ is a projective $A[T]$-module of rank $> \dim A=\dim A[T]-1$,
then $(i)$ $P$ is cancellative and $(ii)$ $P$ has a unimodular
element.

Let $P$ be a projective $A[X_1,\ldots,X_l]$-module of rank $> \dim A$,
then $(i)$ Ravi Rao \cite{Rao} proved that $P$ is cancellative and
$(ii)$ Bhatwadekar and Roy \cite{BRoy} proved that $P$ has a
unimodular element, thus generalizing the Plumstead's result.

Let $P$ be a projective $R=A[X_1,\ldots,X_l,Y_1^{\pm
1},\ldots,Y_m^{\pm 1}]$-module of rank $> \dim A$, then $(i)$
Lindel \cite{lindel} proved that if rank
$P\geq$ max $(2,1+\dim A)$, then $\EL^1(R\oplus P)$ acts transitively on
$\Um(R\oplus P)$. In particular, $P$ is cancellative and $(ii)$
Bhatwadekar, Lindel and Rao \cite{BLR} proved that $P$ has a
unimodular element.

In another direction, Ravi Rao \cite{R2} generalized Plumstead's result by
proving that if $R=A[T,1/g(T)]$ or
$R=A[T,\frac{f_1(T)}{g(T)},\ldots,\frac{f_r(T)}{g(T)}]$, where $g
(T)\in A[T]$ is a non-zerodivisor and if $P$ is a projective
$R$-module of rank $> \dim A$, then $P$ is cancellative. We will generalize
Rao's result by proving that $\EL^1(R\oplus P)$ acts transitively on
$\Um(R\oplus P)$ (\ref{111}). 

 Let $R[X_1,\ldots,X_l,Y_1,\ldots,Y_m, \frac{1}{f_1\ldots f_m}]$,
where $f_i\in A[Y_i]$ and let $P$ be a projective $R$-module of rank
$\geq$ max $\{2,\dim A+1\}$ Then we show that $(i)$ $\EL^1(R\oplus P)$
acts transitively on $\Um(R\oplus P)$ and $(ii)$ If $A$ is an affine
algebra over a field, then $P$ has a unimodular element (\ref{23},
\ref{222}), thus generalizing results of (\cite{lindel}, \cite{BLR})
where it is proved for $f_i=Y_i$.

As an application of the above result, we prove the following result
(\ref{26}): Let $\ol k$ be an algebraically closed field with $1/d!\in
\ol k$ and let $A$ be an affine $\ol k$-algebra of dimension $d$.
Let $R=A[T,1/g(T)]$ or
$R=A[T,\frac{f_1(T)}{g(T)},\ldots,\frac{f_r(T)}{g(T)}]$, where $g(T)$
is a monic polynomial and $g(T),f_1(T),\ldots,f_r(T)$ is
$A[T]$-regular sequence. Then every projective $R$-module of rank
$\geq d$ is cancellative. (See\cite{M2} for motivation)


\section{Preliminaries}  

Let $A$ be a ring and let $M$ be an $A$-module. For $m\in M$, we
define $O_M(m) = \{ \varphi(m) | \varphi \in \Hom_R(M,R) \}$. We say
that $m$ is {\it unimodular} if $O_M(m) = A$. The set of all unimodular
elements of $M$ will be denoted by $\Um(M)$.  We denote by $\Aut_A(M)$,
the group of all $A$-automorphism of $M$.  For an ideal $J$ of $A$, we
denote by $\Aut_A(M,J)$, the kernel of the natural homomorphism
$\Aut_A(M)\rightarrow \Aut_A(M/JM)$.

We denote by $\EL^1(A\oplus M,J)$, the subgroup of $\Aut_A(A\oplus M)$
generated by all the automorphisms $\Delta_{a\varphi}=\left(
\begin{smallmatrix}
 1 & a\varphi\\
0 & id_M
\end{smallmatrix}  \right)
$ and $\Gamma_{m}=\left(\begin{smallmatrix}
1&0\\
m&id_M 
\end{smallmatrix}\right)$
with $a\in J,\varphi \in \Hom_A(M,A)$ and $m \in M$. 

We denote by $\Um^1(A\oplus M, J)$, the set of all $(a,m) \in
\Um(A\oplus M)$ such that $a \in 1 + J$ and by $\Um(A\oplus M,J)$, the
set of all $(a,m) \in \Um^1(A\oplus M,J)$ with $m \in JM$. We will
write $\Um^1_r(A,J)$ for $\Um^1(A\op A^{r-1},J)$ and $\Um_r(A,J)$ for
$\Um(A\op A^{r-1},J)$.

We will write $\EL_{r}^1(A,J)$ for $\EL^1(A\oplus A^{r-1}, J)$,
$\EL^1_r(A)$ for $ \EL_{r}^1(A,A)$ and
$\EL^1(A\oplus M)$ for $\EL^1(A\oplus M,A)$.

\begin{remark}\label{14}
$(i)$ Let $I\subset J$ be ideals of a ring $A$ and let $P$ be a projective
$A$-module. Then, it is easy to see that the natural map $\EL^1(A\oplus
P,J)\rightarrow \EL^1(\frac AI\oplus \frac P{IP},\frac JI)$ is surjective.

$(ii)$ Let $E_r(A)$ be the group generated by elementary matrices
$E_{i_0j_0}(a)=(a_{ij})$, where $i_0\not= j_0$, $a_{i,j}\in A$,
$a_{ii}=1$, $a_{i_0j_0} =a$ and remaining $a_{ij}=0$ for $1\leq
i,j\leq r$. Then using (\cite{W1}, Lemma 2.1), it is easy to see
that $E_r(A)=\EL^1_r(A)$.
\end{remark}

The following result is a consequence of a theorem of Eisenbud-Evans
as stated in (\cite{P}, p.1420).

\begin{theorem}\label{EE}
Let $R$ be a ring and let $P$ be a projective $R$-module of rank
$r$. Let $(a,\ga)\in (R\op P^*)$. Then there exists $\gb\in P^*$ such
that $\hh I_a\geq r$, where $I=(\ga+a\gb)(P)$. In particular, if the
ideal $(\ga(P),a)$ has height $\geq r$, then $\hh I\geq r$. Further,
if $(\ga(P),a)$ is an ideal of height $\geq r$ and $I$ is a proper
ideal of $R$, then $\hh I=r$.
\end{theorem}

The following two results are due to Wiemers (\cite{W1},
Proposition 2.5 and Theorem 3.2).

\begin{proposition}\label{4}
Let $A$ be a ring and let $R = A[X_1,\ldots,X_n,Y_1^{\pm 1},\ldots,Y_m^{\pm
1}]$. Let $c$ be the element $1$, $X_n$ or $Y_m - 1$. If $s \in
A$ and $r\geqq$ max $\{3, dim A + 2\}$, then $\EL_r^1(R,sc)$ acts
transitively on $\Um_r^1(R,sc)$.
\end{proposition}

\begin{theorem}\label{5}
Let $A$ be a ring and let $R = A[X_1,\ldots,X_n,Y_1^{\pm 1},\ldots,Y_m^{\pm
1}]$. Let $P$ be a projective $R$-module of rank $r\geqq$
max $\{2, dim A + 1\}$. If $J$ denotes the ideal $R$, $X_nR$ or $(Y_m
- 1)R$, then $\EL^1(R\oplus P, J)$ acts transitively on $\Um^1(R\oplus P,
J)$.
\end{theorem}

The following result is due to Ravi Rao (\cite{R2}, Lemma 2.1).

\begin{lemma}\label{rr}
Let $B \subset C$ be rings of dimension $d$ and $x \in B$ such that
$B_x = C_x$. Then

$(i)$ $ B/(1 + xb) = C/(1 + xb)$ for all $b \in B$.

$(ii)$ If $I$ is an ideal of $C$ such that $ht I\geq d$ and $I + xC =
C$, then there exists $b \in B$ such that $1 + xb \in I$.  

$(iii)$ If $c \in C$, then $c = 1 + x + x^2 h$ mod $(1 + xb)$ for some
$h \in B$ and for all $b \in B$.
\end{lemma}

\begin{definition}
Let $A$ be a ring and let $M,N$ be $A$-modules. Suppose $f,g :
M\iso N$ be two isomorphisms. We say that ``$f$ is {\it isotopic}
to $g$'' if there exists an isomorphism $\phi: M[X]\iso N[X]$ such
that $\phi(0) = f$ and $\phi(1) = g$. 

Note that if $\sigma\in \EL^1(A\op P)$, then $\sigma$ is isotopic to
identity.
\end{definition}

The following lemma follows from the well known Quillen splitting
lemma (\cite{Q}, Lemma 1) and its proof is essentially contained in
(\cite{Q}, Theorem 1).

\begin{lemma}\label{8}
Let $A$ be a ring and let $P$ be a projective $A$-module. Let $s,t\in
A$ be two comaximal elements. Let $\sigma\in \Aut_{A_{st}}(P_{st})$
which is isotopic to identity. Then $\sigma=\tau_s\theta_t$, where
$\tau\in \Aut_{A_t}(P_t)$ such that $\tau=id$ modulo $sA$ and
$\theta\in \Aut_{A_s}(P_s)$ such that $\theta=id$ modulo $tA$. 
\end{lemma}

The following two results are due to Suslin (\cite{SU}, Corrolary 5.7
and Theorem 6.3).

\begin{theorem}\label{100}
Let $A$ be any ring and let $f\in A[X]$ be a monic polynomial. Let
$\alpha \in \GL_r(A[X])$ be such that $\alpha_f\in
EL^1_r(A[X]_f)$. Then $\alpha\in EL^1_r(A[X])$.
\end{theorem}

\begin{theorem}\label{101}
Let $A$ be a ring and $B=A[X_1,\ldots,X_l]$. Then the canonical map
$\GL_r(B)/EL^1_r(B)\rightarrow K_1(B)$ is an isomorphism for $r\geq$
max $\{3,dimA+2\}$. In particular, if $\ga\in\GL_r(B)$ is stably
elementary, then $\ga$ is elementary.
\end{theorem}


\section{Main Theorem} 

We begin this section with the following result which is easy to
prove. We give the proof for the sake of completeness.

\begin{lemma}\label{19}
Let $A$ be a ring and let $P$ be a projective $A$-module. Let ``bar''
denote reduction modulo the nil radical of $A$.  For an ideal $J$ of
$A$, if $\EL^1(\ol A\oplus \ol P, \ol J)$ acts transitively on $\Um^1(\ol
A\oplus\ol P,\ol J)$, then $\EL^1(A\oplus P,J)$ acts transitively on
$\Um^1(A\oplus P,J)$.
\end{lemma}

\begin{proof} 
Let $(a,p)\in \Um^1( A\oplus P,J)$.  By hypothesis, there exists a
$\sigma\in \EL^1(\ol A\oplus \ol P,\ol J)$ such that $\sigma (\ol a,\ol
p) =(1,0)$. Using (\ref{14}), let $\varphi \in \EL^1(A\oplus P,J)$ be a
lift of $\sigma$ such that $\varphi(a,p)=(1+b,q)$, where $b\in N=
nil(A)$ and $q\in NP$. Note that $b\in N\cap J$. Since $1+b$ is a
unit, we get $\Gamma_1=\Gamma_{\frac{-q}{1+b}} \in \EL^1(A\op P,J)$
such that $\Gamma_1(1+b,q) = (1+b,0)$. It is easy to see that there
exists $p_1,\ldots,p_n \in P$ and $\alpha_1,\ldots, \alpha_n \in P^*$
such that $\alpha_1(p_1)+ \ldots +\alpha_n(p_n) = 1$. Write
$h=\sum_2^n \ga_i(p_i)$.  Note that $(1+b,0)=(1+\sum_1^n
b\ga_i(p_i),0)$, $\Gamma_{\frac{p_1}{1+b}}(1+b,0) = (1+b,p_1)$ and
$\Delta_{-b\alpha_1} (1+b,p_1) =(1+bh,p_1)$, where
$\Delta_{-b\alpha_1}\in \EL^1(A\op P, J)$. Since $1+bh$ is a unit,
$\Gamma_{\frac{-p_1}{1+bh}} (1+bh,p_1)= (1+bh,0)=(1+\sum_2^n
b\ga_i(p_i),0) $. Applying further transformations as above, we can
take $(1+\sum_2^n b\ga_i(p_i),0)$ to $(1,0)$ by elements of $\EL^1(A\op
P,J)$.  $ \hfill \square$
\end{proof}

The following lemma is similar to the Quillen's splitting lemma
(\ref{8}). We will give the sketch of the proof.

\begin{lemma}\label{9} 
Let $A$ be a ring and let $u,v$ be two comaximal elements of $A$.  For
any $s \in A$, every $\alpha \in \EL_{n}^1(A_{uv},s)$ has a splitting
$(\alpha_1)_v \circ(\alpha_2)_u$, where $ \alpha_1\in \EL_{n}^1(A_u,s)$
and $ \alpha_2\in \EL_{n}^1(A_v,s)$.
\end{lemma}

\begin{proof} If $\alpha \in \EL_{n}^1(A_{uv},s)$, then $\alpha = \prod_{i =
1}^r \alpha_i$, where $\alpha_i$ is of the form
$\left(\begin{smallmatrix}
1 & s\underline{v}\\
0 & Id_M \\
\end{smallmatrix}\right)$ or  $\left(\begin{smallmatrix}
1 & 0\\
\underline{w^t} & Id_M \\
\end{smallmatrix}\right),$
where $M=A_{uv}^{n-1}$, $\underline{v},\underline{w}\in M$.

Define $\alpha(X) \in \EL_{n}^1(A[X]_{uv},s)$ by $\alpha(X) =\prod_{i =
1}^r \alpha_i(X)$, where $\alpha_i(X)$ is of the form $
\left(\begin{smallmatrix}
1 & sX\underline{v}\\
0 & Id_{M[X]} \\
\end{smallmatrix}\right)$ or  $\left(\begin{smallmatrix}
1 & 0\\
X\underline{w^t} & Id_{M[X]} \\
\end{smallmatrix}\right)$ as may by the case above.

Since $\alpha(0)=id$ and $\alpha(1)=\ga$, $\alpha$ is isotopic to
identity. Using proof of (\ref{8}) (\cite{M1}, Lemma 2.19), we get
that $\alpha(X) = (\psi_1(X))_v \circ(\psi_2(X))_u$, where $\psi_1(X)
= \alpha(X)\circ\alpha(\lambda u^kX)^{-1}\in \EL_n^1(A_u[X],s)$ and
$\psi_2(X) = \alpha(\lambda u^kX)\in \EL_n^1(A_v[X],s)$ with $\lambda
\in A$, $k \gg 0$. Write $\psi_1(1)=\ga_1 \in \EL_n^1(A_u,s)$ and
$\psi_2(1)=\ga_2 \in \EL^1_n(A_v,s)$, we get that
$\ga(1)=\alpha =
(\ga_1)_v \circ(\ga_2)_u$. 
$\hfill \square$
\end{proof}

\begin{definition}
Let $A$ be a ring of dimension $d$ and let $l,m,n \in
\BN\cup\{0\}$. We say that a ring $R$ is of the type $A\{d,l,m,n\}$,
if $R$ is an $A$-algebra generated by $X_1,\ldots, X_l$, $Y_1,\ldots,
Y_m$, $T_1,\ldots, T_n$, $\frac{1}{f_1\ldots f_m}$,
$\frac{g_{11}}{h_1},\ldots,\frac{g_{1{t_1}}}{h_1},\ldots,\frac{g_{n1}}{h_n},
\ldots,\frac{g_{n{t_n}}}{h_n}$, where $X_i$'s, $Y_i$'s and $T_i$'s are
variables over $A$, $f_i\in A[Y_i]$, $g_{ij}\in A[T_i]$, $h_i\in A[T_i]$
and $h_i$'s are non-zerodivisors.
\end{definition}

For Laurent polynomial ring (i.e. $f_i=Y_i$), the following result is
due to Wiemers (\ref{4}).

\begin{proposition}\label{20}
Let $A$ be a ring of dimension $d$ and let
$R=A[X_1,\ldots,X_l,Y_1,\ldots,Y_m,\frac{1}{f_1\ldots f_m}]$, where
$f_i\in A[Y_i]$ (i.e. $R$ is of the type $A\{d,l,m,0\}$). If $s\in A$ and
$r\geq $ max $\{3, d+2\}$, then $\EL_{r}^1(R,s)$ acts transitively on
$\Um_{r}^1(R,s)$.
\end{proposition}

\begin{proof} 
Without loss of generality, we may assume that $A$ is reduced. The
case $m =0$ is due to Wiemers (\ref{4}).  Assume $m\geq 1$ and
apply induction on $m$.

Let $(a_1,\ldots,a_{r}) \in \Um_{r}^1(R,s)$.  Consider a multiplicative
closed subset $S = 1 + f_mA[Y_m]$ of $A[Y_m]$.  Then $R_S
=B[X_1,\ldots,X_l,Y_1,\ldots,Y_{m-1},\frac{1}{f_1\ldots f_{m-1}}]$, where
$B=A[Y_m]_{{f_{m}S}}$ and $dim B = dim A$. Since $R_S$
is of the type $B\{d,l,m-1,0\}$, by induction hypothesis on $m$, there
exists  $\sigma \in \EL_{r}^1(R_S,s)$ such that
$\sigma (a_1,\ldots,a_{r})= (1,0,\ldots,0)$. We can find $g\in S$ and
$\sigma' \in \EL_{r}^1(R_g,s)$ such that $\sigma' (a_1,\ldots,a_{r})=
(1,0,\ldots,0)$.

Write
$C=A[X_1,\ldots,X_l,Y_1,\ldots,Y_m,\frac{1}{f_1,\ldots,f_{m-1}}]$. Consider
the following fiber product diagram

\[
\xymatrix{
C \ar@{->}[r]
     \ar@{->}[d]
& R=C_{f_m} 
     \ar@{->}[d] 
\\
C_g   
\ar@{->}[r]
     &R_g=C_{gf_m}.     
}
\]

Since $\sigma'\in \EL_{r}^1(C_{gf_m},s)$, by (\ref{9}), $\sigma'
=(\sigma_2)_{f_m}\circ (\sigma_1)_{g} $, where $\sigma_2 \in
\EL_{r}^1(C_g,s )$ and $\sigma_1 \in \EL_{r}^1(R,s)$. Since
$(\sigma_1)_{g}(a_1,\ldots,a_{r})=(\sigma_2)_{f_m}^{-1}
(1,0,\ldots,0)$, patching $\sigma_1(a_1,\ldots,a_{r}) \in
\Um^1_r(C_{f_m},s)$ and $(\sigma_2)^{-1} (1,0,\ldots,0) \in
\Um^1_r(C_g,s)$, we get a unimodular row $(c_1,\ldots,c_{r}) \in
\Um_{r}^1(C,s)$.  Since $C$ is of the type $A\{d,l+1,m-1,0\}$, by
induction hypothesis on $m$, there exists $\phi \in \EL_{r}^1(C,s)$
such that $\phi (c_1,\ldots,c_{r}) = (1,0,\ldots,0)$. Taking
projection, we get $\Phi \in \EL_{r}^1(R,s)$ such that
$\Phi\sigma_1(a_1,\ldots,a_{r}) = (1,0,\ldots,0)$. This completes the
proof.  $ \hfill \square$
\end{proof}

\begin{corollary}\label{60}
Let $A$ be a ring of dimension $d$ and let
$R=A[X_1,\ldots,X_l,Y_1,\ldots,Y_m,\frac{1}{f_1\ldots f_m}]$, where
$f_i\in A[Y_i]$.  Let $c$ be $1$ or $X_l$, If $s\in A$ and $r\geq $
max $\{3, d+2\}$, then $\EL_{r}^1(R,sc)$ acts transitively on
$\Um_{r}^1(R,sc)$.
\end{corollary}

\begin{proof}
Let $(a_1,\ldots,a_r)\in \Um_{r}^1(R,sc)$. The case $c=1$ is done by
(\ref{20}).  Assume $c=X_l$. We can assume, after an
$\EL_r^1(R,sX_l)$-transformation, that $a_2,\ldots,a_r \in
sX_lR$. Then we can find $(b_1,\ldots,b_r)\in
\Um_r(R,sX_l)$ such that the following equation holds:
$$a_1b_1+\ldots +a_rb_r = 1. \hspace*{1in} (i)$$
Now consider the $A$-automorphism $\mu : R\rightarrow R$ defined as follows
\begin{center}
$ X_i\mapsto X_i$ for $i= 1,...,l-1$,
\\$X_l \mapsto X_l(f_1\ldots f_m)^N$ for some large positive integer $N$.
                           
\end{center}
Applying $\mu$, we can read the image of equation $(i)$ in the subring
$S=A[X_1,\ldots,X_l,Y_1,\ldots,Y_m]$. By (\ref{4}), we obtain $\psi\in
\EL_r^1(R,sX_l)$ such that $\psi(\mu(a_1),\ldots,\mu(a_r)) =
(1,0,\ldots,0)$.  Since $\mu^{-1}(X_l)$ and $X_l$ generate the same
ideal in $R$, applying $\mu^{-1}$, the proof follows.
 $\hfill \square$
\end{proof}

\begin{corollary}\label{200}
Let $A$ be a ring of dimension $d$ and let
$R=A[X_1,\ldots,X_l,Y_1,\ldots,Y_m,\frac{1}{f_1\ldots f_m}]$, where
$f_i\in A[Y_i]$. Then $\EL^1_r(R)$ acts transitively on $\Um_r(R)$ for
$r\geq max\{3,d+2\}$.
\end{corollary}

The following result is similar to (\cite{R2}, Theorem 5.1).

\begin{theorem}\label{201} 
Let $A$ be a ring of dimension $d\geq 1$ and let
$R=A[X_1,\ldots,X_l,Y_1,\ldots,Y_m,\frac{1}{f_1\ldots f_m}]$, where
$f_i\in A[Y_i]$ (i.e. $R$ is of the type $A\{d,l,m,0\}$). Then 

(i) the canonical map $\Phi_r: GL_r(R)/\EL^1_r(R)\rightarrow K_1(R)$ is
surjective for $r\geq d+1$.

(ii) Assume $f_i\in A[Y_i]$ is a monic polynomial for all $i$. Then for
$r \geq $ max $\{3, d+2\}$, any stably elementary matrix in $GL_r(R)$ is
in $\EL^1_r(R)$. In particular, $\Phi_{d+2}$ is an isomorphism.
\end{theorem}

\begin{proof} $(i)$ Let $[M] \in K_1(R)$. We have to show that 
$[M]=[N]$ in $K_1(R)$ for some $N \in GL_{d+1}(R)$. Without loss of
generality, we may assume that $M \in
GL_{d+2}(R)$. By (\ref{20}), there exists an elementary matrix $\sigma
\in \EL^1_{d+2}(R)$ such that $ M\sigma =
 \begin{pmatrix}
 M' & a\\
0 & 1
\end{pmatrix}$. 
Applying further $\sigma'\in \EL^1_{d+2}(R) $, we get $\sigma'M\sigma =
 \begin{pmatrix} N & 0\\ 0 & 1
\end{pmatrix}$, where $M', N \in GL_{d+1}(R)$. Hence $[M]=[N]$ in $K_1(R)$. 
This completes the proof of $(i)$.

$(ii)$ Let $M \in GL_r(R)$ be a stably elementary matrix.  For $m=0$,
we are done by (\ref{101}). Assume $m\geq 1$.

Let $S = 1 + f_mA[Y_m]$.  Then $R_S =
B[X_1,\ldots,X_l,Y_1,\ldots,Y_{m-1},\frac{1}{{f_1\ldots f_{m-1}}}]$,
where $B= A[Y_m]_{{f_{m}S}}$ and $dim B = dim A$.  
Since $R_S$ is of the type $B\{d,l,m-1,0\}$, by induction
hypothesis on $m$, $M \in \EL^1_r(R_S)$. Hence there exists $g\in S$ such
that $M\in \EL^1_r(R_g)$. Let $\sigma \in
\EL^1_r(R_g)$ be such that $\sigma M = Id$. 

Write $C=A[X_1,\ldots,X_l,Y_1,\ldots,Y_m,\frac{1}{f_1\ldots
f_{m-1}}]$.  Consider the following fiber product diagram
$$ \xymatrix{ C
\ar@{->}[r] \ar@{->}[d] & C_{f_m}=R \ar@{->}[d] \\
C_g \ar@{->}[r]
&C_{gf_m}=R_g.}$$ 
By (\ref{9}), $\sigma = (\sigma_2)_{f_m}\circ (\sigma_1)_g$,
where $\sigma_2 \in \EL^1_{r}(C_g)$ and $\sigma_1 \in \EL^1_{r}(C_{f_m})$. 
Since
$(\sigma_1 M)_{g} = (\sigma_2)_{f_m}^{-1}$, patching $\sigma_1M$ and
$(\sigma_2)^{-1}$, we get $N\in GL_{r}(C)$ such that $N_{f_m}=\sigma_1 M$. 

Write $D=A[X_1,\ldots,X_n,Y_1,\ldots,Y_{m-1},\frac{1}{f_1\ldots
f_{m-1}}]$. Then $D[Y_m]=C$ and $D[Y_m]_{f_m}=R$. Since $N\in
GL_r(D[Y_m])$, $f_m\in D[Y_m]$ is a monic polynomial and
$N_{f_m}=\sigma_1 M$ is stably elementary, by 
(\ref{100}), $N$ is stably elementary.  Since $C$ is of
the type $A\{d,l+1,m-1,0\}$, by induction hypothesis on $m$, $N\in
\EL^1_r(C)$. Since $\sigma_1$ is elementary, we get that $M \in
\EL^1_r(R)$. This completes the proof of $(ii)$.  $ \hfill \square$
\end{proof}

\begin{lemma}\label{22}
Let $R$ be a ring of the type $A\{d,l,m,n\}$. Let $P$ be a projective
$R$-module of rank $r \geq$ max $\{2, 1+d\}$. Then there exists an
$s\in A$, $p_1,\ldots ,p_r \in P$ and $\varphi_1,\ldots ,\varphi_r \in
\Hom(P,R)$ such that the following properties holds.

(i) $P_s$ is free.  

(ii) $(\varphi_i(p_j) ) =$ diagonal $(s,s,\ldots ,s)$.

(iii) $sP \subset p_1A + \ldots  +p_rA$.

(iv) The image of $s$ in $A_{red}$ is a nonzero divisor.  

(v) $(0: sA) = (0: s^2A)$.
\end{lemma}

\begin{proof} 
Without loss of generality, we may assume that $A$ is reduced. Let $S$
be the set of all non-zerodivisors in $A$. Since $dim A_S=0$ and
projective $R_S$-module $P_S$ has a constant rank, we may assume that $A_S$
is a field. Then it is easy to see that
$A_S[T_i,\frac{g_{ij}}{h_i}]=A_S[T_i,\frac{1}{h_i}]$ (assuming gcd
$(g_{ij}, h_i)=1$).  Therefore $R_S=
A_S[X_1,\ldots,X_l,Y_1,\ldots,Y_m,T_1,\ldots ,T_n, \frac{1}{f_1\ldots
f_mh_1\ldots h_n}]$ is a localization of a polynomial
ring over a field. Hence projective modules over $R_S$ are stably free.
Since $P_S$ is stably free of rank $\geq$ max $\{2,1+d\}$, by
(\ref{20}), $P_S$ is a free $R_S$-module of rank $r$. We can
find an $s \in S$ such that $P_s$ is a free $R_s$-module. The
remaining properties can be checked by taking a basis $p_1,\ldots ,p_r
\in P$ of $P_s$, a basis $\varphi_1,\ldots ,\varphi_r \in \Hom(P,R)$ of
$P_s^*$ and by replacing $s$ by some power of $s$, if needed. This
completes the proof. $ \hfill \square$
\end{proof}

\begin{lemma}\label{15}
Let $R$ be a ring of the type $A\{d,l,m,n\}$. Let $P$ be a projective
$R$-module of rank $r$. Choose $s\in A$, $p_1,\ldots ,p_r \in P$ and
$\varphi_1,\ldots ,\varphi_r \in \Hom(P,R)$ satisfying the properties
of (\ref{22}). Let $(a,p)\in \Um(R\oplus P,sA)$ with $p = c_1p_1
+\ldots + c_rp_r$, where $c_i\in sR$ for $i=1$ to $r$. Assume there
exists $\phi \in \EL_{r+1}^1(R,s)$ such that $\phi(a,c_1,\ldots
,c_r) = (1,0,\ldots ,0)$. Then there exists $ \varPhi \in \EL^1(R\oplus P)$
such that $\varPhi (a,p) = (1,0)$.
\end{lemma}

\begin{proof} Since $\phi \in \EL_{r+1}^1(R,s)$, $\phi
=\prod_{j=1}^{n}\phi_j$, where $\phi_j = \Delta_{s\psi_j}$ or
$\Gamma_{v^t}$ with $\psi_j = (b_{1j},\ldots ,b_{rj})\in {R^r}^*$ and $v
=(f_1,\ldots, f_r)\in {R^{r}} $. 

Define a map $\Theta :
\EL_{r+1}^1(R,s) \rightarrow \EL^1(R\oplus P)$ as follows
$$ \Theta (\Delta_{s\psi_j}) =
\begin{pmatrix}
 1 & \sum_{i=1}^r {b_{ij}}\varphi_i\\
0 & id_P
\end{pmatrix}  \;\;\; and \;\;\;\Theta (\Gamma_{v^t}) =\begin{pmatrix}
1&0\\
\sum_{i=1}^r f_ip_i &id_P 
\end{pmatrix}.$$ 

Let $\varPhi = \prod_{j=1}^{n}\Theta(\phi_j) \in \EL^1(R\op P)$. Then
it is easy to see that $\varPhi(a,p) = (1,0)$. This completes the
proof.  $ \hfill \square$
\end{proof}

\begin{remark}\label{15.1}
From the proof of above lemma, it is clear that if
$\phi \in \EL_{r+1}^1(R,sX_l)$ such that $\phi(a,c_1,\ldots ,c_r) =
(1,0,\ldots ,0)$, then $\Phi \in \EL^1(R\oplus P,X_l)$ such that
$\varPhi (a,p) = (1,0)$.  
\end{remark}

For Laurent polynomial ring (i.e. $f_i=Y_i$ and $J=R$), the following
result is due to Lindel \cite{lindel}.

\begin{theorem}\label{23}
Let $A$ be a ring of dimension $d$ and let
$R=A[X_1,\ldots,X_l,Y_1,\ldots,Y_m,\frac{1}{f_1\ldots f_m}]$, where
$f_i\in A[Y_i]$ (i.e. $R$ is of the type $A\{d,l,m,0\}$). Let $P$ be a
projective $R$-module of rank $r \geq $ max $\{2, d+1\}$. If $J$
denote the ideal $R$ or $X_lR$, then $\EL^1(R\oplus P,J)$ acts
transitively on $\Um^1(R\oplus P,J)$.
\end{theorem}

\begin{proof} 
Without loss of generality, we may assume that $A$ is a reduced. We
will use induction on $d$. When $d =0$, we may assume that $A$ is a
field. Hence projective modules over $R$ are stably free (proof
of lemma \ref{22}). Using (\ref{60}), we are done.

Assume $d >0$. By (\ref{22}), there exists a non-zerodivisor
$s\in A$, $p_1,\ldots ,p_{r} \in P$ and $\phi_1,\ldots ,\phi_{r} \in
P^* = \Hom_R(P,R)$ satisfying the properties of (\ref{22}). If
$s\in A$ is a unit, then $P$ is a free and the result follows from
(\ref{60}). Assume $s$ is not a unit.

Let $(a,p) \in \Um^1(R\oplus P,J)$. Let ``bar'' denotes reduction modulo
the ideal $s^2R$. Since $dim\ol A< dim A$, by induction hypothesis,
there exists $\varphi \in \EL^1(\ol R\oplus \ol P,J)$ such that
$\varphi(\ol a,\ol p) = (1,0)$. Using (\ref{14}), let $\Phi \in
\EL^1(R\oplus P,J)$ be a lift of $\varphi$ and $\Phi(a,p)= (b,q)$,
where $ b\equiv 1$ mod $s^2JR$ and $q \in s^2JP$.

By (\ref{22}), there
exists $a_1,\ldots ,a_{r} \in sJR$ such that $q = a_1p_1 + \ldots +
a_{r}p_{r}$. It follows that $(b,a_1,\ldots ,a_{r}) \in
\Um_{r+1}(R,sJ)$. By (\ref{60}), there exists
$\phi \in \EL_{r+1}^1(R,sJ)$ such that $\phi(b,a_1,\ldots
,a_{r})=(1,0,\ldots ,0)$. Applying (\ref{15.1}), we get $\Psi \in
\EL^1(R\oplus P,J)$ such that $\Psi(b,q) = (1,0)$. Therefore
$\Psi\Phi(a,p)=(1,0)$. This completes the proof.  $ \hfill \square$
\end{proof}

For Laurent polynomial ring (i.e. $f_i=Y_i$), the following result is
due to Bhatwadekar-Lindel-Rao \cite{BLR}. 

\begin{theorem}\label{222}
Let $k$ be a field and let $A$ be an affine $k$-algebra of dimension
$d$. Let $R=A[X_1,\ldots,X_l,Y_1,\ldots,Y_m,\frac{1}{f_1\ldots f_m}]$,
where $f_i\in A[Y_i]$ (i.e. $R$ is of the type $A\{d,l,m,0\}$). Then
every projective $R$-module $P$ of rank $\geq d+1$ has a unimodular
element.
\end{theorem}

\begin{proof} 
We assume that $A$ is reduced and use induction on $\dim
A$. If $\dim A=0$, then every projective module of constant rank is
free (\ref{20}, \ref{22}). Assume $\dim A>0$.  

By (\ref{22}), there exists a non-zerodivisor $s\in A$ such that $P_s$
is free $R_s$-module. Let ``bar'' denote reduction modulo the ideal
$sR$. By induction hypothesis, $\ol P$ has a unimodular element, say $\ol
p$. Clearly $(p,s)\in \Um(P\oplus R)$, where $p\in P$ is a lift of $\ol
p$. By (\ref{EE}), we may assume that $\hh I \geq d+1$, where
$I=O_P(p)$. We claim that $I_{(1+sA)} = R_{(1+sA)}$ (i.e. $p\in
\Um(P_{1+sA})$).

Since $R$ is a Jacobson ring, $\sqrt I=\cap \mm$ is the intersection
of all maximal ideals of $R$ containing $I$. Since $I+sR=R$, $s \notin
(I\cap A)$. Let $\mm$ be any maximal ideal of $R$ which contains
$I$. Since $A$ and $R$ are affine $k$-algebras, $\mm\cap A$ is a maximal
ideal of $A$. Hence $\mm \cap A$ contains an element of the form
$1+sa$ for some $a\in A$ (as $s\notin \mm\cap A$). Hence $\mm
R_{(1+sA)}=R_{(1+sA)}$ and $I_{(1+sA)}=R_{(1+sA)}$. This proves the claim.

Let $S= 1+sA$. Choose $t\in S$ such that $p\in \Um(P_t)$.
Let $p_1\in \Um(P_s)$ and $p\in \Um(P_t)$. Since $R_{sS}$ is of the type
$A_{sS}\{d-1,l,m,0\}$, by (\ref{23}), there exist $\varphi \in
\EL^1(P_{sS})$ such that $\varphi(p_1) = p$. We can choose $t_1=tt_2\in
S$ such that $\varphi \in \EL^1(P_{st_1})$.  By (\ref{8}), $\varphi =
(\varphi_1)_s\circ(\varphi_2)_{t_1}$, where $\varphi_2\in \EL^1(P_s)$
and $\varphi_1 \in \EL^1(P_{t_1})$.
Consider the following fiber product diagram
$$
\xymatrix{
P \ar@{->}[r]
     \ar@{->}[d]
&P_s 
     \ar@{->}[d]
\\
P_{t_1} \ar@{->}[r]
     &P_{st_1}.     
}$$

Since $(\varphi_2)_{t_1}(p_1) = (\varphi_1)_s^{-1}(p)$, patching
$\varphi_2(p_1) \in \Um (P_s)$ and $\varphi_1^{-1}(p) \in \Um(P_t)$,
we get a unimodular element in $P$.  This proves the result. $ \hfill
\square$
\end{proof}

The following result generalizes a result of Ravi Rao \cite{R2} where it is
proved that $P$ is cancellative.

\begin{theorem}\label{111}
Let $A$ be a ring of dimension $d$ and let
$R=A[X,\frac{f_1}g,\ldots,\frac{f_n}g]$, where $g,f_i\in A[X]$ with
$g$ a non-zerodivisor. Let $P$ be a projective $R$-module of rank
$r\geq$ max $\{2,d+1\}$. Then $\EL^1(R\op P)$ acts transitively on
$\Um(R\op P)$.
\end{theorem}

\begin{proof}
We will assume that $A$ is reduced and apply induction on $\dim A$. If
$\dim A=0$, then we may assume that $A$ is a field. Hence $R$ is a PID
and $P$ is free. By (\ref{4}), we are done.

Assume $\dim A=d>0$. By (\ref{22}), we can choose a non-zerodivisor
$s\in A$, $p_1,\ldots,p_r\in P$ and $\phi_1,\ldots,\phi_r\in P^*$
satisfying the properties of (\ref{22}).

Let $(a,p)\in \Um(R\op P)$.  Let ``bar'' denotes reduction modulo
$sgR$. Then $\dim \ol R <\dim R$ and $r\geq \dim \ol R +1$. By Serre's
result \cite{Serre}, $\ol P$ has a unimodular element, say $\ol
q$. Then $(0,\ol q)\in \Um(\ol R\op \ol P)$. By Bass result \cite{Bas},
there exists $\phi\in \EL^1(\ol R\op \ol P)$ such that $\phi(\ol a,\ol
p)=(0,\ol q)$. Using (\ref{14}), let $\Phi\in \EL^1(R\op P)$ be a lift
of $\phi$ and $\Phi(a,p)=(b,q)$, where $b\in sgR$. By (\ref{EE}), we
may assume that $\hh O_P(q) \geq d+1$.

Write $B=A[X], x=sg,I=O_P(q)$ and $C=R$. Then $\dim B=\dim C$ and
$B_{sg}=C_{sg}$. By (\ref{rr}(ii)), there exists $h\in A[X]$ such that
$1+sgh\in O_P(q)$. Hence there exists $\gf\in P^*$ such that
$\gf(q)=1+sgh$.

By (\ref{rr}(iii)), there exists $b'\in R$ such that
$b-b'(1+sgh)=1+sg+s^2g^2h'$ for some $h'\in A[X]$. Since
$\Delta_{-b'\gf}(b,q)=(b-b'\gf(q),q)=(1+sg+s^2g^2h',q)=(b_0,q)$ and
$\Gamma_{-q}(b_0,q)=(b_0,q-b_0q)=(b_0,sgq_1)$ for some $q_1\in P$ and
$b_0\in A[X]$ with $b_0=1$ mod $sgA[X]$.

Write $sgq_1=c_1p_1+\ldots+c_rp_r$ for some $c_i\in R$. Then
$(b_0,c_1,\ldots,c_r)\in \Um^1_{r+1}(R,sg)$. It is easy to see that by
adding some multiples of $b_0$ to $c_1,\ldots,c_r$, we may assume that
$(b_0,c_1,\ldots,c_r)\in \Um^1(A[X],sgA[X])$.  By (\ref{4}), there
exists $\Theta\in \EL^1_{r+1}(A[X],s)$ such that $\Theta(b_0,c_1,\ldots,
c_r)=(1,0,\ldots,0)$.
Applying (\ref{15}), there exists $\Psi\in \EL^1(R\op P)$ such that
$\Psi(b_0,sgq_1)=(1,0)$. This proves the result.
$\hfill \square$.
\end{proof}

\begin{question}
Let $R$ be a ring of type $A\{d,l,m,n\}$ and let $P$ be a projective
$R$-module of rank $\geq $max $\{2,d+1\}$. 

$(i)$ Does $\EL^1(R\op P)$ acts transitively on $\Um(R\op P)$? In
particular, Is $P$ cancellative?

$(ii)$ Does $P$ has a unimodular element?
\end{question}  

Assume $n=0$. Then $(i)$ is (\ref{23}) and for affine algebras over a
field, $(ii)$ is (\ref{222}).

When either $P$ is free or $\ol k=\ol \BF_p$, then the following
result is proved in \cite{M2}.

\begin{theorem}\label{26}
Let $\ol k$ be an algebraically closed field with $1/d! \in \ol k$ and
let $A$ be an affine $\ol k$-algebra of dimension $d$.  Let
$f(T) \in A[T]$ be a monic polynomial and assume that either

$(i)$ $ R = A[T,\frac{1}{f(T)}]$ or 

$(ii)$ $R= A[T,\frac{f_1}{f},\ldots ,\frac{f_n}{f}]$, where
$f,f_1,\ldots ,f_n$ is $A[T]$-regular sequence. \\
Then every projective $R$-module $P$ of rank $d$ is cancellative.
\end{theorem}

\begin{proof} 
By (\ref{22}), there exists a non-zerodivisor $s\in A$ satisfying the
properties of (\ref{22}).  Let $(a,p) \in \Um(R\oplus P)$.

Let ``bar'' denote reduction modulo ideal $s^3A$. Since $dim \ol A<
dim A$, by (\ref{23}, \ref{111}), there exists a $\phi\in \EL^1(\ol
R\oplus \ol P)$ such that $\phi (\ol a,\ol p) = (1,0)$. Let $\Phi \in
\EL^1(R\oplus P)$ be a lift of $\phi$. Then $\Phi(a,p) = (b,q)$, where
$(b,q) \in \Um^1(R\oplus P, s^2A)$. Now the proof follows by
(\cite{M2}, Theorem 4.4).  $ \hfill \square$
\end{proof}

The proof of the following result is same as of (\ref{26}) using
(\cite{M2}, Theorem 5.5).

\begin{theorem}\label{270}
Let $k$ be a real closed field and let $A$ be an affine $k$-algebra of
dimension $d-2$. Let $f\in A[X,T]$ be
a monic polynomial in $T$ which does not belong to any real maximal
ideal of $A[X,T]$. Assume that either 

$(i)$ $R= A[X,T,1/f]$ or 

$(ii)$ $R = A[X,T,f_1/f,\ldots ,f_n/f]$, where $f,f_1,\ldots ,f_n$ is
$A[X,T]$-regular sequence.\\ Then every projective $R$-module of rank
$d-1$ is cancellative.
\end{theorem}

\section{An analogue of Wiemers result}

We begin this section with the following result which can be proved
using the same arguments as in (\cite{W1}, Corollary 3.4) and using
(\ref{23})

\begin{theorem}\label{24}
Let $A$ be a ring of dimension $d$ and $R = A[X_1,\ldots
,X_l,Y_1,\ldots ,Y_m,\frac{1}{f_1\ldots f_m}]$, where $f_i\in
A[Y_i]$. Let $P$ be a projective $R$-module of rank $\geq d+1$. Then
the natural map $\Aut_R(P)\rightarrow \Aut_{\ol R}(P/X_lP)$ with $\ol R
= R/X_lR$ is surjective.
\end{theorem}

Using the automorphism $\mu$ defined in (\ref{60}), the following
result can be proved using the same arguments as in (\cite{W1},
Proposition 4.1).

\begin{proposition}\label{25}
Let $A$ be a ring of dimension $d$, $1/d! \in A$ and $R = A[X_1,\ldots
,X_l,Y_1,\ldots ,Y_m,\frac{1}{f_1\ldots f_m}]$ with $l\geq 1$, $f_i\in
A[Y_i]$. Then $GL_{d+1}(R,X_lJR)$ acts transitively
on $\Um_{d+1}(R,X_lJR)$, where $J$ is an ideal of $A$.
\end{proposition}

When $f_i=Y_i$, the following result is due to Wiemers (\cite{W1},
Theorem 4.3).  The proof of this result is same as of (\cite{W1},
Theorem 4.3) using (\ref{24}, \ref{25}).

\begin{theorem}\label{261}
Let $A$ be a ring of dimension $d$ with $1/d!\in A$ and let $R =
A[X_1,\ldots ,X_l,Y_1,\ldots ,Y_m,\frac{1}{f_1\ldots f_m}]$ with
$f_i\in A[Y_i]$ for $i =1$ to $m$ . Let $P$ be a projective $R$-module
of rank $\geq d$. If $Q$ is another projective $R$-module such that
$R\oplus P \cong R\oplus Q$ and $\ol P \cong \ol Q$, then $P\cong Q$,
where ``bar'' denote reduction modulo the ideal $(X_1,\ldots ,X_l)R$ .
\end{theorem}

Using (\ref{26}, \ref{261}), we get the following result.

\begin{corollary}
Let $\ol k$ be an algebraically closed field with $1/d! \in \ol k$ and
let $A$ be an affine $\ol k$-algebra of dimension $d$.  Let $f(T) \in
A[T]$ be a monic polynomial and let $R=A[X_1,\ldots,X_l,T,\frac
1{f(T)}]$. Then every projective $R$-module of rank $\geq d$ is
cancellative.
\end{corollary}

{}

\end{document}